\documentclass[12pt]{amsart}

\usepackage{amssymb}
\usepackage{latexsym}
\usepackage{amsmath}

\begin{document}

\def\c{{\mathbb C}}
\def\n{{\mathbb N}}

\newtheorem{theo}{Theorem}

\title[Quantum groups and Fuss-Catalan algebras]{Quantum groups and Fuss-Catalan algebras}

\author{Teodor Banica}

\address{Institut de Math{\'e}matiques de Jussieu, 175 rue du Chevaleret, 75013 Paris}

\email{banica@math.jussieu.fr}

\maketitle

\centerline{Institut de Math{\'e}matiques de Jussieu, 175 rue du
  Chevaleret, 75013 Paris, France}
\centerline{banica@math.jussieu.fr}

\bigskip

{\bf Abstract.} The categories of representations of compact quantum groups of
automorphisms of certain inclusions of finite dimensional
$\c^*$-algebras are shown to be isomorphic to the categories of
Fuss-Catalan diagrams.

\section*{Introduction}

Let $(D,\tau )$ be a finite dimensional $\c^*$-algebra together
with a trace. In \cite{w2} Wang constructs an algebra
$A_{aut}(D,\tau )$: the biggest Hopf $\c^*$-algebra
co-acting on $(D,\tau )$.

From the point of view of noncommutative geometry the algebra $D$
corresponds to a noncommutative finite space
$\widehat{D}$ and the trace $\tau$ corresponds to a mesure
$\widehat{\tau}$ on $\widehat{D}$. Thus $A_{aut}(D,\tau )$
corresponds to the ``quantum symmetry group''
$G_{aut}(\widehat{D},\widehat{\tau})$.

If $D=\c^n$ with $n=1,2,3$ then $A_{aut}(D,\tau )$ is just the
algebra of functions on the $n$-th symmetric group. But if $dim(D)\geq
4$ then $A_{aut}(D,\tau )$ is infinite-dimensional (\cite{w2}). 

The corepresentations of $A_{aut}(D,\tau )$ are studied in
\cite{aut}: under a suitable assumption on the trace $\tau$, the algebras of
symmetries of the fundamental corepresentation (i.e. the one on $D$)
are shown to be isomorphic to the Temperley-Lieb algebras.

Let $B\subset D$ be an inclusion of finite dimensional $\c^*$-algebras
and let $\varphi$ be a state on $D$. Following Wang one can construct
a Hopf algebra $A_{aut}(B\subset D,\varphi )$: the biggest Hopf
$\c^*$-algebra co-acting on $D$ such that $B$ and $\varphi$ are left
invariant.

The main result in this paper is that, under suitable assumptions on
$\varphi$, the category of finite dimensional corepresentations of
$A_{aut}(B\subset D,\varphi )$ is isomorphic to the completion of the
category of Fuss-Catalan diagrams. (These are certain colored Temperley-Lieb
diagrams, discovered by Bisch and Jones in connection with intermediate
subfactors \cite{bj}.) The proof (\S 1 -- \S 4) uses \cite{aut},
\cite{bj}, \cite{wo2} and reconstruction techniques.

The Fuss-Catalan diagrams were recently shown to appear in several
contexts, related to subfactors, planar algebras and integrable
lattice models. See e.g. Bisch and Jones \cite{bj}, \cite{bj2}, Di
Francesco \cite{df}, Landau \cite{la} and the references there in. In
the last section of the paper (\S 5) we discuss the relation between
$A_{aut}(B\subset D,\varphi )$ and subfactors.

\section{Preliminaries}

The Fuss-Catalan category, as well as other categories to be used in
this paper, is a tensor $\c^*$-category having $(\n ,+)$ as monoid
of objects. For simplifying writing, such a tensor category
will be called a $\n$-algebra. If $C$ is a $\n$-algebra we use the notations
$${C} (m,n)=Hom_C (m,n)\hskip 2cm {C}(m)=End_C(m)$$

As a first class of examples, associated to any object $O$ in a tensor
$\c^*$-category is the $\n$-algebra $\n O$ given by
$$\n O(m,n)=Hom(O^{\otimes m}, O^{\otimes n})$$

Fix $\delta >0$. The $\n$-algebra $TL^2$ is defined as follows. The space $TL^2(m,n)$ consists of linear combinations
of Temperley-Lieb diagrams between $2m$ points and $2n$ points
$${{TL^2}}(m,n)=\left\{ \sum \,\alpha\,\,\begin{matrix}\cdots\cdots & \leftarrow &
    2m\,\,\mbox{points}\cr \mathfrak{W} & \leftarrow &
    m+n\mbox{ strings}\cr \cdot \,\cdot & \leftarrow &
    2n\,\,\mbox{points}\end{matrix}\right\}$$
(strings join points and don't cross) and the operations
$\circ$, $\otimes$ and $*$ are induced by vertical and horizontal
concatenation and upside-down turning of diagrams. With the following
rule: erasing a circle is the same as multiplying by $\delta$.
$$\mbox{A}\circ \mbox{B}=\begin{matrix}\mbox{B}\cr \mbox{A}\end{matrix}\hskip 2cm \mbox{A}\otimes \mbox{B}
=\mbox{A}\mbox{B}\hskip 2cm \mbox{A}^*=\forall\hskip 2cm\bigcirc
=\delta$$

Consider the following two elements ${u}\in {{TL^2}} (0,1)$ and ${m}\in {{TL^2}} (2,1)$:
$${u} =\delta^{-\frac{1}{2}}\,\,\cap\hskip 2cm {m}
=\delta^{\frac{1}{2}}\,\,\mid\cup\mid$$ 

\begin{theo}
The following relations

(i) ${m}{m}^*=\delta^2 1$

(ii) ${u}^*{u} =1$

(iii) ${m} ({m}\otimes 1)={m} (1\otimes{m} )$

(iv) ${m} (1\otimes{u} )={m} ({u}\otimes 1)=1$

(v) $({m}\otimes 1)(1\otimes{m}^*)=(1\otimes{m} )({m}^*\otimes 1)={m}^*{m}$

\noindent are a presentation of ${{TL^2}}$ by ${u}\in {{TL^2}} (0,1)$ and ${m}\in {{TL^2}} (2,1)$.
\end{theo}

The result says that $u$ and $m$ satisfy the relations and that if
${C}$ is a $\n$-algebra and $v\in C(0,1)$ and $n\in C(2,1)$ satisfy
the relations then there exists a $\n$-algebra morphism
$$TL^2\to{C}\hskip 2cm u\mapsto v\hskip 2cm m\mapsto n$$

This is proved in \cite{aut}. Actually in \cite{aut} the ``index''
$\delta^2$ is an integer and $u$ and $m$ are certain explicit
operators, but these extra structures are not used.

Let $D$ be a finite dimensional $\c^*$-algebra with a state $\varphi$
on it. We have a scalar product $<x,y>=\varphi (y^*x)$ on $D$, so $D$
is an object in the category of finite dimensional Hilbert
spaces. Consider the unit $u$ and the multiplication $m$ of $D$.
$${u}\in\n D (0,1)\hskip 2cm {m}\in \n D (2,1)$$

The relations in theorem 1 are satisfied if and only if the first one,
namely ${m}{m}^*=\delta^2 1$, is satisfied. If $D=\oplus M_{n_\beta}$ is a
decomposition of $D$, we must have $Tr(Q_\beta^{-1})=\delta^2$ for
any block $Q_\beta$ of the unique $Q\in D$ such that $\varphi
=Tr(Q.)$. This can be cheked by direct computation; see \cite{aut} for the case
$\varphi$ = trace.

A linear form $\varphi$ such that ${m}{m}^*=\delta^2 1$, where the
adjoint of the multiplication is taken with respect to the scalar product
associated to $\varphi$, will be called a $\delta$-form.

One can deduce from theorem 1 that if $\varphi$ is a $\delta$-form
then category of corepresentations of the Hopf $\c^*$-algebra $A_{aut}(D,\varphi )$ is the
completion of ${TL^2}$. The case $\varphi$ = trace is studied in
\cite{aut}; for the general case, see \S 4 below.

\section{The Fuss-Catalan category}

A Fuss-Catalan diagram is a planar diagram formed by an upper row of
$4m$ points, a lower row of $4n$ points and by $2m+2n$
non-crossing strings joining them. Both rows of points are
colored from left to right in the following standard way
$$\mbox{white, black, black, white, white, black, black,...}$$
and strings have to join pairs of points having the same color.

Fix $\beta >0$ and $\omega >0$. The $\n$-algebra $FC$ is defined as
follows. The spaces $FC(m,n)$ consist of linear combinations of Fuss-Catalan diagrams
$${{FC}}(m,n)=\left\{ \sum \,\alpha\,\,\begin{matrix}\mbox{\tiny{w,b,b,w,w,b,b,w,\ldots}} & \leftarrow &
    4m\,\,\mbox{colored points}\cr \ & \ &
    m+n\mbox{ white strings}\cr \mathfrak{W}\mbox{\ \ } & \leftarrow & \mbox{and}\cr  \ & \ &
    m+n\mbox{ black strings}\cr \mbox{\tiny{w,b,b,w,w,b,b,w,\ldots}} & \leftarrow &
    4n\,\,\mbox{colored points}\end{matrix}\right\}$$
and the operations
$\circ$, $\otimes$ and $*$ are induced by vertical and horizontal
concatenation and upside-down turning of diagrams. With the following
rule: erasing a black/white circle is the same as
multiplying by $\beta$/$\omega$.
$$\mbox{A}\circ \mbox{B}=\begin{matrix}\mbox{B}\cr \mbox{A}\end{matrix}\hskip 10mm \mbox{A}\otimes \mbox{B}
=\mbox{A}\mbox{B}\hskip 10mm \mbox{A}^*=\forall\hskip 10mm\mbox{\tiny{black}}\rightarrow\bigcirc =\beta\hskip 10mm
\mbox{\tiny{white}}\rightarrow\bigcirc =\omega$$

Let $\delta =\beta\omega$. The following bicolored analogues of the elements $u$ and $m$ in \S 1
$${u} =\delta^{-\frac{1}{2}}\,\,
\bigcap\;\!\!\!\!\!\!\!\!\cap\hskip 2cm {m} =\delta^{\frac{1}{2}}\,\,
\mid\mid\bigcup\;\!\!\!\!\!\!\!\!\cup\mid\mid$$
generate in ${FC}$ a $\n$-subalgera which is isomorphic to
$TL^2$.

Consider also the black and white Jones projections.
$${e}
=\omega^{-1}\,\,\mid\begin{matrix}\cup\cr\cap\end{matrix}\mid\hskip 2cm
{f}
=\beta^{-1}\,\,\mid\mid\mid\begin{matrix}\cup\cr\cap\end{matrix}\mid\mid\mid$$

We have ${f} =\beta^{-2}(1\otimes{m}{e} ){m}^*$, so we won't need $f$ for presenting ${FC}$.

For simplifying writing we identify  ${x}$ and ${x}\otimes 1$, for any ${x}$.

\begin{theo}
The following relations, with ${f} =\beta^{-2}(1\otimes{m}{e} ){m}^*$
and $\delta =\beta\omega$

(1) the relations in theorem 1

(2) ${e} ={e}^2={e}^*$, ${f} ={f}^*$ and $(1\otimes {f} ){f} ={f} (1\otimes {f} )$

(3) ${e}{u} ={u}$

(4) ${m}{e}{m}^*={m} (1\otimes{e} ){m}^*=\beta^21$

(5) ${m}{m} ({e}\otimes{e}\otimes{e} )={e}{m}{m} ({e}\otimes 1\otimes{e} )$

\noindent are a presentation of ${{{{FC}}}}$ by ${m}\in {{{{FC}}}} (2,1)$, ${u}\in {{{{FC}}}} (0,1)$ and ${e}\in {{{{FC}}}} (1)$.
\end{theo}

\begin{proof}
As for any presentation result, we have to prove two assertions.

(I) The elements $m,u,e$ satisfy the relations (1--5) and generate the
$\n$-algebra $FC$.

(II) If  $M$, $U$ and $E$ in a $\n$-algebra $C$ satisfy the relations
(1--5), then there exists a morphism of $\n$-algebras $FC\to C$
sending ${m}\mapsto M$, ${u}\mapsto U$ and ${e}\mapsto E$.

The proof will be based on the results in the
paper of Bisch and Jones \cite{bj}, plus some graphic computations for (I)
and some purely algebraic computations for (II).

(I) First, the relations (1--5) are easily verified by drawing pictures.

Let us show that the $\n$-subalgebra
${C}=<m,u,e>$ of $FC$ is equal to ${FC}$. First, $C$ contains the
infinite sequence of black and white Jones projections
$${p}_1={e}=\omega^{-1}\,\,\mid\begin{matrix}\cup\cr\cap\end{matrix}\mid$$
$${p}_2={f}=\beta^{-1}\,\,\mid\mid\mid\begin{matrix}\cup\cr\cap\end{matrix}\mid\mid\mid$$
$${p}_3=1\otimes{e} =\omega^{-1}\,\,\mid\mid\mid\mid\mid\begin{matrix}\cup\cr\cap\end{matrix}\mid$$
$$p_4=1\otimes f
=\beta^{-1}\,\,\mid\mid\mid\mid\mid\mid\mid\begin{matrix}\cup\cr\cap\end{matrix}\mid\mid\mid$$
$$\ldots$$
as well as the infinite sequence of bicolored Jones projections
$${e}_1={u}{u}^*=\delta^{-1}\,\,\begin{matrix}\bigcup\hskip -3.9mm\cup\cr
  \bigcap\hskip -3.9mm\cap\end{matrix}$$
$$e_2=\delta^{-2}{m}^*{m}=\delta^{-1}\,\,\mid\mid\begin{matrix}\bigcup\hskip -3.86mm\cup\cr\bigcap\hskip -3.86mm\cap\end{matrix}\mid\mid$$
$${e}_3=1\otimes{u}{u}^*=\delta^{-1}\,\,\mid\mid\mid\mid\begin{matrix}\bigcup\hskip -3.86mm\cup\cr
  \bigcap\hskip -3.86mm\cap\end{matrix}$$
$${e}_4=\delta^{-2}(1\otimes
{m}^*{m}
)=\delta^{-1}\,\,\mid\mid\mid\mid\mid\mid\begin{matrix}\bigcup\hskip -3.86mm\cup\cr\bigcap\hskip -3.86mm\cap\end{matrix}\mid\mid$$
$$\ldots$$
which by the result of Bisch and Jones (\cite{bj}) generate the
diagonal $\n$-algebra $\Delta {{{{FC}}}}$. (If $X$ is a $\n$-algebra,
its diagonal $\Delta X$ is defined by $\Delta X(m)=X(m)$ and $\Delta 
X(m,n)=\emptyset$ if $m\neq n$.) Thus we have inclusions
$$\Delta {{{{FC}}}}\subset{C}\subset{{{{FC}}}}$$
so we can use the following standard argument. First, we have $\Delta
{{{{FC}}}} =\Delta C$. Second, the existence of semicircles shows that
the objects of $C$ and ${FC}$ are selfdual and by Frobenius reciprocity
we get
$$dim({C} (m,n))=dim\left( {C} \left( \frac{m+n}{2}\right) \right)
=dim\left( {{{{FC}}}} \left( \frac{m+n}{2}\right) \right)=dim({{{{FC}}}}
(m,n))$$
for $m+n$ even. By tensoring with  $u$ and $u^*$ we get embeddings
$$C(m,n)\subset C(m,n+1)\hskip 2cm {{FC}}(m,n)\subset {{FC}}(m,n+1)$$
and this shows that the dimension equalities hold for any $m$ and
$n$. Together with $\Delta {{{{FC}}}}\subset{C}\subset{{{{FC}}}}$, this
shows that ${C} ={{{{FC}}}}$.

(II) Assume that $M$, $U$ and $E$ in a $\n$-algebra $C$ satisfy the
relations (1--5). We have to construct a morphism ${{{{FC}}}}\to C$ sending
$${m}\mapsto M ,\, {u}\mapsto U ,\, {e}\mapsto E$$

This will be done in two steps. First, we restrict attention to
diagonals: we would like to construct a morphism $\Delta
{{{{FC}}}}\to\Delta C$ sending
$$m^*m\mapsto M^*M,\,
uu^*\mapsto UU^* ,\, {e}\mapsto E$$

By constructing the corresponding Jones projections $E_i$ and $P_i$, we must send
$${e}_i\mapsto E_i ,\, {p}_i\mapsto P_i \,\,\,\,\, (i=1,2,3,\ldots
)$$

The presentation result for $\Delta {{{{FC}}}}$ of Bisch and Jones
(\cite{bj}) reduces this to an algebraic computation. More precisely,
it is proved in \cite{bj} that the following relations

\smallskip

(a) ${e}_i^2 = {e}_i$, ${e}_i {e}_j = {e}_j {e}_i$ if $|i-j| \ge 2$ and ${e}_i {e}_{i \pm 1} {e}_i = \delta^{-2}{e}_i$

(b) ${p}_i^2 = {p}_i$ and ${p}_i {p}_j = {p}_j {p}_i$

(c) ${e}_i{p}_i = {p}_i {e}_i = {e}_i$ and ${p}_i {e}_j = {e}_j {p}_i$ if $|i-j| \ge 2$

(d) ${e}_{2i \pm 1}{p}_{2i}{e}_{2i \pm 1} = \beta^{-2} {e}_{2i \pm 1}$ and ${e}_{2i}{p}_{2i \pm 1}{e}_{2i} = \omega^{-2} {e}_{2i}$

(e) ${p}_{2i}{e}_{2i \pm 1}{p}_{2i} = \beta^{-2} {p}_{2i \pm 1}{p}_{2i}$ and ${p}_{2i \pm 
    1}{e}_{2i} {p}_{2i \pm 1} = \omega^{-2} {p}_{2i} {p}_{2i \pm 1}$

\smallskip

\noindent are a presentation of $\Delta {{{{FC}}}}$. So it remains to
verify that
$$(1-5)\Longrightarrow (a-e)$$
where ${m}$, ${u}$ are ${e}$ are abstract objects and we are no
longer allowed to draw pictures.

First, by using ${e}_{n+2}=1\otimes {e}_n$ and ${p}_{n+2}=1\otimes {p}_n$
these relations reduce to:

\smallskip

($\alpha$) ${e}_i^2={e}_i$ for $i=1,2$, $e_1e_2e_1=\delta^{-2}e_1$ and
$e_2e_1e_2=\delta^{-2}e_2$.

($\beta$) ${p}_i^2={p}_i$ for $i=1,2$ and $[{p}_1,{p}_2]=[1\otimes {p}_1,{p}_2]=[1\otimes {p}_2,{p}_2]=0$

($\gamma$) $[{e}_2,1\otimes {p}_2]=[{p}_2,1\otimes {e}_2]=0$ and ${e}_i{p}_i={p}_i{e}_i={e}_i$ for $i=1,2$

($\delta 1$) ${e}_1{p}_2{e}_1=\beta^{-2}{e}_1$ and $(1\otimes
{e}_1){p}_2(1\otimes {e}_1)=\beta^{-2}(1\otimes {e}_1)$

($\delta 2$) ${e}_2{p}_1{e}_2={e}_2(1\otimes {p}_1){e}_2=\omega^{-2}{e}_2$

($\epsilon 1$) $\beta^2{p}_2{e}_1{p}_2=\omega^2{p}_1{e}_2{p}_1=
{p}_1{p}_2$

($\epsilon 2$) $\beta^2{p}_2(1\otimes {e}_1){p}_2=\omega^2(1\otimes
{p}_1){e}_2(1\otimes {p}_1)=(1\otimes {p}_1){p}_2$

\smallskip

With ${e}_1={u}{u}^*$, ${e}_2=\delta^{-2}{m}^*{m}$, ${p}_1={e}$
and ${p}_2={f}$ one can see that most of them are trivial. What is
left can be reformulated in the following way.

\smallskip

(x) ${e}{m}^*{m} {e} =\beta^2{f}^*{e}$

(y) $(1\otimes {e} ){m}^*{m} (1\otimes {e} )=\beta^2{f}^*(1\otimes {e} )$

(z) ${f}^*={f}^*{f}$

(t) $[{e} ,{f} ]=[1\otimes {e} ,{f} ]=[{m}^*{m} ,1\otimes
{f} ]=[{f} ,1\otimes {m}^*{m} ]=0$

\smallskip

By multiplying the relation (5) by ${u}$ and by $1\otimes 1\otimes{u}$ to the right
we get the following useful formula, to be used many times in what follows.
$${m} ({e}\otimes {e} )={e}{m} (1\otimes {e} )={e}{m}{e}$$

Let us verify (x--t). First, we have
$$\beta^2{f}^*{e} ={m} ({e}\otimes {e} )(1\otimes{m}^*)$$
and by replacing ${m} ({e}\otimes {e} )$ with ${e}{m}{e}$ we get
${e}{m}^*{m}{e}$, so (x) is true. We have
$$(1\otimes {e} ){m}^*{m} (1\otimes {e} )={m} (1\otimes ({e}{m} (1\otimes {e}
))^*)$$
and by replacing ${e}{m} (1\otimes {e} )$ with ${e}{m}{e}$ we get
$\beta^2{f}^*(1\otimes {e} )$, so (y) is true. We have
$${f}^*{f} =\beta^{-4}{m} (1\otimes {e}{m}^*{m} {e} ){m}^*$$
and by replacing ${e}{m}^*{m} {e}$ with ${e}{m} {e} (1\otimes{m}^* )$,
then ${e}{m} {e} $ with ${m} ({e}\otimes {e} )$ we get
${f}^*$, so (z) is true. The first two commutators are zero, because ${f}{e}$ and
${f}(1\otimes {e} )$ are selfadjoint. Same for the others, because of
the formulae
$${m}{m}^*(1\otimes {f}{f}^*)=\beta^{-4}(1\otimes
1\otimes {m} {e} ){m}^*{m}^*{m}{m} (1\otimes 1\otimes {e}{m}^*)$$ 
$$(1\otimes{m}^*{m}){f}{f}^*=\beta^{-4}(1\otimes {m}^*{m}{e} ){m}^*{m}(1\otimes
{e}{m}^*{m})$$

The conclusion is that we constructed a certain $\n$-algebra morphism
$$\Delta J:\Delta {{{{FC}}}}\to\Delta C$$
that we have to extend now to a morphism $J:{{FC}}\to C$ sending $u\mapsto U$ and
$m\mapsto M$. We will use a standard argument (see \cite{kw}). For $w$
bigger than $k$ and $l$ we define
$$\phi :{{FC}}(l,k)\rightarrow {{FC}}(w)\hskip 1cm x\mapsto 
(u^{\otimes (w-k)}\otimes 1_k)\, x\, ((u^*)^{\otimes (w-l)}\otimes 1_l)$$ 
$$\theta :{{FC}}(w)\rightarrow {{FC}}(l,k)\hskip 18mm x\mapsto
((u^*)^{\otimes (w-k)}\otimes 1_k)\, x\, (u^{\otimes (w-l)}\otimes 1_l)$$
where $1_k=1^{\otimes k}$ and where the convention $x=x\otimes 1$ is
no longer used. We define $\Phi$ and
$\Theta$ in $C$ by similar formulae. We have $\theta\phi =\Theta\Phi
=Id$. We define a map $J$ by
$$\begin{matrix}{{FC}}(l,k)& \displaystyle{\mathop{\longrightarrow}^{J}}&
  C(l,k)\cr \phi\downarrow & \ & \uparrow \Theta\cr
{{FC}}(w)& \displaystyle{\mathop{\longrightarrow}^{\Delta J}}&
  C(w)\end{matrix}$$

As $J(a)$ doesn't depend on the choice of $w$, these $J$'s are the
components of a map $J:{{FC}}\to C$. This map $J$ extends $\Delta J$ and sends
$u\mapsto U$ and $m\mapsto M$. It remains to prove that $J$ is a
morphism. We have
$$Im(\phi )=\{ x\in {{FC}}(w)\mid 
x=((uu^*)^{\otimes (w-k)}\otimes 1_k)\, x\, ((uu^*)^{\otimes
  (w-l)}\otimes 1_l)\}$$
as well as a similar description of $Im(\Phi )$, so $J$ sends $Im(\phi
)$ to $Im(\Phi )$. On the other hand we have $\Theta\Phi =Id$, so $\Phi\Theta =Id$ on $Im(\Phi )$. Thus
$$\begin{matrix}{{FC}}(l,k)& \displaystyle{\mathop{\longrightarrow}^{J}}&
  C(l,k)\cr \phi\downarrow & \ & \downarrow \Phi\cr
{{FC}}(w)& \displaystyle{\mathop{\longrightarrow}^{\Delta J}}&
  C(w)\end{matrix}$$
commutes, so $J$ is multiplicative:
$$J(ab)=\Theta (\Delta J\phi (a)\Delta J\phi (b))=\Theta (\Phi
J(a)\Phi J(b))=\Theta\Phi (J(a)J(b))=J(a)J(b)$$

It remains to prove that $J(a\otimes b)=J(a)\otimes J(b)$. We have
$a\otimes b=(a\otimes 1_s)(1_t\otimes b)$ for certain $s$ and $t$, so
it is enough to prove it for pairs $(a,b)$ of the form $(1_t,b)$ or
$(a,1_s)$. For $(a,1_s)$ this is clear, so it remains to prove that the set
$$B=\left\{ b\in {{FC}}\mid J(1_t\otimes b)=1_t\otimes J(b),\,\forall\,
  t\in\n\right\}$$
is equal ${{FC}}$. First, $\Delta J$ being a $\n$-algebra morphism, we have $\Delta
   {{FC}}\subset B$. On the other hand, computation gives $J(1_t\otimes u \otimes 1_s)=1_t\otimes U \otimes
   1_s$. Also, $J$ being involutive and multiplicative, $B$ is stable
   by involution and multiplication. We conclude that $B$ contains the
   compositions of elements of $\Delta
   {{FC}}$ with $1_t\otimes u\otimes 1_s$'s and $1_t\otimes u^*
   \otimes 1_s$'s. But any $b$ in ${FC}$ is equal to $\theta\phi (b)$,
   so it is of this form and we are done.
\end{proof}

\section{Inclusions of finite dimensional $\c^*$-algebras}

Let $B\subset D$ be an inclusion of finite dimensional
$\c^*$-algebras and let $\varphi$ be a state on $D$. We have the scalar
product $<x,y>=\varphi (y^*x)$ on $D$. The multiplication $m$ of $D$, the
unit $u$ of $D$ and the orthogonal projection $e$ from $D$ onto $B$
$${m} :D\otimes D\to D\hskip 1cm {u}
:{\c{}}\to D\hskip 1cm {e} :D\to D$$
can be regarded as elements of the $\n$-algebra $\n D$ given by $\n D(m,n)={\mathcal
  L}(D^{\otimes m}, D^{\otimes n})$.

We say that $\varphi$ is a $(\beta ,\omega )$-form on $B\subset D$ if
it is a $\beta\omega$-form on $D$, if its restriction $\varphi_{\mid B}$ is a
$\beta$-form on $B$ and if ${e}$ is a $B-B$ bimodule morphism. (For
$\delta$-forms, see \S 1.)

As a first example, if $\phi$ is a $\beta$-form on $B$ and $\psi$ is a $\omega$-form on $W$
then $\phi\otimes\psi$ is a $(\beta ,\omega )$-form on $B\subset
B\otimes W$. In particular a $\delta$-form on $D$ is a $(1,\delta
)$-form on $\c\subset D$.

\begin{theo}
If $\varphi$ is a $(\beta ,\omega )$-form on $B\subset D$ then
$<m,u,e>= {FC}$.
\end{theo}

\begin{proof}
We prove that $m,u,e$ satisfy the relations (1--5). The formulae ${e}={e}^2={e}^*$ and (3) are true, (1) is equivalent to
the fact that $\varphi$ is a $\beta\omega$-form and (5) says that
$${e} (b){e} (c){e} (d)={e} ({e} (b)c{e} (d))$$
for any $b,c,d$ in $B$, i.e. that ${e}$ is a morphism of $B-B$ bimodules.

Let $\{ b_{-i}\}_{i\in\n}$ be an orthonormal basis of $B$ and let
$\{b_j\}_{j\in\n}$ be an orthonormal basis of $B^\perp$. We denote by
$\{ b_n\}_{n\in{{\mathbb Z}}}$ the orthonormal basis $\{
b_{-i},b_j\}_{i,j\in\n}$ of $D$. We have
$${m}^*(b)=\sum_{k,s\in{{\mathbb Z}}}b_k\otimes b_s<b,b_kb_s>=\sum_{k,s\in{{\mathbb Z}}}b_k\otimes
b_s<b_k^*b,b_s>=\sum_{k\in{{\mathbb Z}}}b_k\otimes b_k^*b$$
so $\varphi$ is a $\delta$-form if and only if
$\sum b_kb_k^*=\delta^21$. On the other hand, we get
$${m} {e}{m}^*(b)={m}\left( \sum_{k\in{{\mathbb Z}}}{e} (b_k)\otimes b_k^*b\right) =\left(
  \sum_{k\in{{\mathbb Z}}} {e} (b_k)b_k^*\right) b=\left(
  \sum_{i\in\n} b_{-i}b_{-i}^*\right) b$$
so the formula ${m} {e}{m}^*=\beta^21$ in (4) is equivalent to the
fact that $\varphi_{\mid N}$ is a $\beta$-form on $B$.

It remains to check the following three formulae, with ${f}
=\beta^{-2}(1\otimes{m}{e} ){m}^*$. 
$${f} ={f}^*\hskip 1cm (1\otimes {f} ){f} ={f} (1\otimes{f} )\hskip
1cm {m} (1\otimes {e} ){m}^*=\beta^21\hskip 1cm (\star)$$

By using the fact that $e$ is a bimodule morphism we get succesively that
$$\sigma (B)=B\hskip 2cm e*=*e$$
where $\sigma :D\to D$ is such that $\varphi (ab)=\varphi (b\sigma
(a))$. By using the above formula for ${m}^*$ we get
$${f} (x\otimes y)=\beta^{-2}(1\otimes{m} {e} ){m}^*(x\otimes
y)=\beta^{-2}\sum_{k\in{{\mathbb Z}}}b_k\otimes {e} (b_k^*b_m)b_n$$

This allows us to prove the first $(\star )$ formula, because we have
$$<{f} (b_m\otimes b_n),b_M\otimes b_N>=\beta^{-2}\varphi
(b_N^*{e} (b_M^*b_m)b_n)=<b_m\otimes b_n,{f} (b_M\otimes b_N)>$$
for any $m,n,M,N$. The second $(\star )$ formula follows from
$$<(1\otimes {f} ){f} (x\otimes y\otimes z),b_k\otimes b_s\otimes
w>=\beta^{-4}<{e} (b_s^*ay)z,w>$$
$$<f(1\otimes {f} )(x\otimes y\otimes z),b_k\otimes b_s\otimes
w>=\beta^{-4}\sum_{t\in{{\mathbb Z}}} <ab_t,b_s><e(b_t^*y)z,w>$$
with $a=e(b_k^*x)$, for any $x,y,z,w,k,s$. For the third $(\star )$
formula, we have
$${m}^*(b)=\sum_{k,s\in{{\mathbb Z}}}b_k\otimes
b_s<b,b_kb_s>=\sum_{k,s\in{{\mathbb Z}}}b_k<b\sigma (b_s^*),b_k>\otimes b_s=
\sum_{s\in{{\mathbb Z}}}b\sigma (b_s^*)\otimes b_s$$
and this gives ${m} (1\otimes {e} ){m}^* (b)=bq$ with $q$ given by
$$q=\sum_{s\in{{\mathbb Z}}} \sigma (b_s^*) {e}
(b_s)=\sum_{i\in{{\mathbb N}}} \sigma (b_{-i}^*) b_{-i}={m}_B{m}_B^*(1)=\beta^21$$
where ${m}_B$ is the multiplication of $N$, which was computed in a
similar way.

Thus $m,u,e$ satisfy the relations (1--5), so theorem 2 applies and gives
a certain $\n$-algebra surjective morphism $J:{{FC}}\to\, <m,u,e>$.

It remains to prove that $J$ is faithful. For, consider the maps
$$\phi_n:FC(n)\to FC(n-1)\hskip 2cm x\mapsto (1^{\otimes (n-1)}\otimes
v^*)(x\otimes 1)(1^{\otimes (n-1)}\otimes v)$$
$$\psi_n:C(n)\to C(n-1)\hskip 14mm x\mapsto (1^{\otimes (n-1)}\otimes
J(v)^*)(x\otimes 1)(1^{\otimes (n-1)}\otimes J(v))$$
where $v=m^*u\in FC(0,2)$. These make the following diagram
commutative
$$\begin{matrix}FC(n)& \displaystyle{\mathop{\longrightarrow}^J}&
  C(n)\cr \phi_n\downarrow & \ & \downarrow \psi_n\cr
FC(n-1)& \displaystyle{\mathop{\longrightarrow}^J}&
  C(n-1)\end{matrix}$$
and by gluing such diagrams we get a factorisation by $J$ of the
composition on the left of conditional expectations, which is the
Markov trace. By positivity $J$ is faithful on $\Delta FC$,
then by Frobenius reciprocity faithfulness has to hold on the
whole $FC$.
\end{proof}

\section{Quantum automorphism groups of inclusions}

Let $B\subset D$ be an inclusion of finite dimensional
$\c^*$-algebras and let $\varphi$ be a state on $D$. Following Wang (\cite{w2}) we define the universal ${\c{}}^*$-algebra
$A_{aut}(B\subset D,\varphi )$ generated by the coefficients $v_{ij}$
of a unitary matrix $v$ subject to the conditions
$${m} \in Hom(v^{\otimes 2},v)\hskip 1cm {u}\in Hom(1,v)\hskip
1cm{e}\in End(v)$$
where $m:D\otimes D\to D$ is the multiplication, ${u}
:{\c{}}\to D$ is the unit and ${e} :D\to D$ is the projection onto $B$, with respect to the scalar product $<x,y>=\varphi (y^*x)$.

This definition has to be understood as follows.  Let $n=dim(D)$ and
fix a vector space isomorphism $D\simeq \c^n$. Let $F$ be the free
$*$-algebra on $n^2$ variables $\{ v_{ij}\}_{i,j=1,...,n}$ and let
$v=(v_{ij})\in M_n\otimes F$. For any $k\in\n$ define $v^{\otimes k}$
to be
$$v^{\otimes k}=v_{1,k+1}v_{2,k+1}\ldots v_{k,k+1}\in M_n^{\otimes
  k}\otimes F$$
If $a,b\in n$ and $t\in {\mathcal L}(M_n^{\otimes a}, M_n^{\otimes
  b})$, the collection of relations between $v_{ij}$'s and their adjoints
obtained by identifying coefficients in the formula 
$$(t\otimes id)v^{\otimes a}=v^{\otimes b}(t\otimes id)$$
can be called ``the relation $t\in Hom(v^{\otimes a}, v^{\otimes
  b})$''. With this definition, let $J\subset F$ be the two-sided $*$-ideal generated by the relations
${m} \in Hom(v^{\otimes 2},v)$, ${u}\in Hom(1,v)$ and ${e}\in End(v)$,
together with the relations obtained by identifying coefficients in
$$vv^*=v^*v=1$$
The $*$-algebra $F/J$ being generated by the coefficients of a unitary
matrix, its enveloping ${\bf C}^*$-algebra of $F/J$ is
well-defined. We call it $A_{aut}(B\subset D,\varphi )$.

By universality we can construct a $\c^*$-algebra morphism
$$\Delta :A_{aut}(B\subset D,\varphi )\to A_{aut}(B\subset D,\varphi
)\otimes A_{aut}(B\subset D,\varphi )$$
sending $v_{ij}\mapsto\sum_k v_{ik}\otimes v_{kj}$ for any $i$ and
$j$ (the tensor product being the ``min'' tensor product). We have
$(id\otimes \Delta )v=v_{12}v_{13}$, so the comultiplication relation
$$(id\otimes\Delta )\Delta
=(\Delta\otimes id)\Delta$$
holds on the generating set $\{
v_{ij}\}_{i,j=1,...,n}$. It follows that the comultiplication relation
holds on the whole $A_{aut}(B\subset
D,\varphi )$. Summing up, we have constructed a pair
$$(A_{aut}(B\subset D,\varphi ),v)$$
consisting of a unital Hopf
$\c^*$-algebra together with a generating corepresentation, i.e. a
compact matrix pseudogroup in the sense of Woronowicz (\cite{wo1}).

The matrix $v$ is a corepresentation
of $A_{aut}(B\subset D,\varphi )$ on the Hilbert space $D$. The three
``Hom'' conditions translate into the fact that $v$ corresponds to a
coaction of $A_{aut}(B\subset D,\varphi )$ on the
${\c{}}^*$-algebra $D$, which leaves $\varphi$ and $B$ invariant. See
Wang (\cite{w2}) and \cite{aut} for details and comments, in the case
$B=\c$.

See \S 3 in \cite{survey} for a general construction of such universal Hopf
$\c^*$-algebras.

We recall from Woronowicz (\cite{wo2}) that the completion of a tensor
$\c^*$-category is by definition the smallest semisimple tensor
$\c^*$-category containing it.

\begin{theo}
If $\varphi$ is a $(\beta ,\omega )$-form on $B\subset D$ then the tensor
$\c^*$-category of finite dimensional corepresentations of $A_{aut}(B\subset
D,\varphi )$ is the completion of ${{{{FC}}}}$.
\end{theo}

\begin{proof}
The unital Hopf $\c^*$-algebra $A_{aut}(B\subset D,\varphi )$ being
presented by the relations corresponding to $m$, $u$ and
$e$, its tensor $\c^*$-category of corepresentations has to be
completion of the tensor $\c^*$-category $<m,u,e>$ generated by $m$,
$u$ and $e$. (This is a direct consequence of tannakian duality
\cite{wo2}, cf. theorem 3.1 in \cite{survey}.)

On the other hand, the linear form $\varphi$ being a $(\beta ,\omega
)$-form, theorem 3 applies and gives an isomorphism $<m,u,e>\simeq FC$.
\end{proof}

In the case $B={\c{}}$ and $\varphi$ = trace, studied in \cite{aut}, we have $FC={TL^2}$.

Note the following corollary of theorem 4: if $\varphi$ is a $(\beta
,\omega )$-form then $\n v={FC}$. (This follows from the basic facts
about completion in \cite{wo2}.)
 
Note also, as an even weaker version of theorem 4, the dimension equalities
$$dim(Hom(v^{\otimes m},v^{\otimes n}))=dim({FC}(m,n))$$
for any $m$ and $n$. Together with standard reconstruction tricks (see
e.g. \cite{survey}) and with the results in \cite{bj}, these
equalities could be used for classifying the irreducible
corepresentations of $A_{aut}(B\subset D,\varphi )$ and for finding
their fusion rules.

\section{Fixed point subfactors}

If $D$ is a finite dimensional $\c^*$-algebra then there
exists a unique $\delta$-trace on it: the canonical trace $\tau_D$. This is by
definition the restriction to $D$ of the unique unital trace of
matrices, via the left regular representation. We have $\delta =\sqrt{dim(D)}$. See \cite{aut}.

In \cite{kac} we construct inclusions of fixed point von Neumann
algebras of the form
$$(P\otimes B)^K\subset (P\otimes D)^K$$
where $P$ is a $II_1$-factor, $B\subset D$ is an inclusion
of finite dimensional ${\c{}}^*$-algebras with a trace $\tau$ and $K$ is
a compact quantum group acting minimally on $P$ and acting on $D$ such
that $B$ and $\tau$ are left invariant. (See \cite{kac} for technical
details, in terms of unital Hopf $\c^*$-algebras.) We also show that
if $K$ acts on $(C,\tau )$ and acts minimally on a $II_1$-factor $P$
then we have the following implications
$$Z((P\otimes C)^K)=\c\Longleftrightarrow Z(C)^K=\c\Longrightarrow
\tau =\tau_C$$
and we deduce from this that we have the following sequence of implications
$$\left\{ \begin{matrix}(P\otimes B)^K\subset (P\otimes D)^K\cr
    \mbox{is a subfactor}\end{matrix}\right\} \Longleftrightarrow
 \left\{ \begin{matrix}Z(D)^K=\c\cr
    Z(B)^K=\c\end{matrix} \right\} \Longrightarrow \left\{
  \begin{matrix}\tau =\tau_D\cr
    \tau_{\mid B}=\tau_B\end{matrix} \right\}$$
which can be glued to the following sequence of elementary implications
$$\left\{ \begin{matrix}\tau =\tau_D\cr
    \tau_{\mid B}=\tau_B\end{matrix} \right\} \Longrightarrow
\left\{ {\tau_D}_{\mid B}=\tau_B  \right\} \Longrightarrow\left\{ Ind(B\subset D)=\frac{dim(D)}{dim(B)}\in\n \right\}$$

The following question is raised in \cite{kac}: is ${\tau_D}_{\mid B}=\tau_B$ the only
restriction on $B\subset D$?

\begin{theo}
For an inclusion $B\subset D$ the following are equivalent

-- there exist subfactors of the form $(P\otimes B)^K\subset (P\otimes D)^K$

-- $B\subset D$ commutes with the canonical traces of $B$ and $D$.
\end{theo}

\begin{proof}
The canonical trace $\tau_D$ is a $\delta$-form, with
$\delta=\sqrt{dim(D)}$. Its restriction to $B$ is the canonical
trace of $B$, so it is a $\beta$-form, with $\beta
=\sqrt{dim(B)}$. Also $\tau_D$ being a trace, the projection $e$ has to
be a $B-B$ bimodule morphism. Thus $\tau_D$ is a $(\beta ,\omega )$-form, with
$$\beta =\sqrt{dim(B)}\hskip 2cm \omega =\sqrt{dim(D)/dim(B)}$$

Thus theorem 4 applies to $\tau_D$. In terms of quantum groups, we get
that the fundamental representation $\pi$ of $K=G_{aut}(B\subset
D,\tau_D )$ satisfies
$$dim(Hom(\pi^{\otimes m},\pi^{\otimes n}))=dim({FC}(m,n))$$
for any $m$ and $n$. With $m=0$ and $n=1$ we get
$$dim(Hom(1,\pi ))=dim({FC}(0,1))=1$$

Together with the canonical isomorphism $D^K\simeq Hom(1,\pi )$, this
shows that the action of $K$ on $D$ is ergodic. In particular 
$Z(D)^K=Z(B)^K=\c$, so by the above, if $P$ is a $II_1$-factor with
minimal action of $K$ (the existence of such a $P$ is shown by Ueda in
\cite{ueda}) then $(P\otimes B)^K\subset (P\otimes D)^K$ is a
subfactor and we are done.
\end{proof}

Note that by \cite{kac} the standard invariant of $P^K\subset (P\otimes
D)^K$ is the Popa system associated to $\pi$, which by theorem 4 is the Fuss-Catalan Popa
system. Equivalently,
$$P^K\subset (P\otimes B)^K\subset (P\otimes D)^K$$
is isomorphic to a free composition of $A_\infty$ subfactors.

In \cite{subf} and \cite{kac} standard invariant of subfactors are
associated to actions of compact quantum groups on objects of the form
$(\c\subset M_n,\varphi )$ and $(B\subset D, \tau )$. This gives
evidence for the existence of a general construction, 
starting with objects of the form $(B\subset D,\varphi
)$, subject to certain conditions. The results in this paper suggest
that there should be only one condition on $(B\subset D,\varphi
)$, namely ``$\varphi$ is a $(\beta ,\omega )$-form on $B\subset D$''.

\end{document}